\documentclass[twoside]{article}
\def\nkyear{200*}
\def\nkvol{**}
\def\nkno{*}

\def\firstpage{1}
\def\lastpage{19}
\def\shorttitle{CURVATURE ESTIMATES FOR IRREDUCIBLE SYMMETRIC SPACES}
\def\shortauthor{{\it Liu Xusheng}}

\usepackage{amsmath}
\usepackage{amsfonts}
\usepackage{amssymb}
\usepackage{amscd}
\usepackage{amsthm}
\usepackage{amsbsy}
\usepackage{graphicx}
\catcode`@=11
\def\ps@nk{\def\@oddhead{\vbox{\hbox to \hsize{\footnotesize \rm \hfill \shorttitle
\hfill \thepage} \vspace{1mm} \hrule \vspace*{-2mm}}}
\def\@evenhead{\vbox{\hbox to \hsize{\footnotesize \rm \thepage \hfill \shortauthor 
\hfill} \vspace{1mm} \hrule \vspace*{-2mm}}}
\def\@oddfoot{} \def\@evenfoot{}}
\def\ps@first{\def\@oddhead{\vbox{\vspace*{-5mm}
\hbox to \hsize{\footnotesize \it Chin. Ann. Math. \hfill} \vspace{1mm}
\hbox to \hsize{\footnotesize \rm {\bf \nkvol B}:\nkno(\nkyear),\thepage--\lastpage. \hfill} \break}}
\def\@evenhead{\vbox{\vspace*{-5mm}
\hbox to \hsize{\footnotesize \it Chin. Ann. Math. \hfill} \vspace{1mm}
\hbox to \hsize{\footnotesize \rm {\bf \nkvol B}:\nkno(\nkyear),\thepage--\lastpage. \hfill} \break}}
\def\@oddfoot{} \def\@evenfoot{}}
\def\Section#1{\Sec{\large #1} \setcounter{equation}{0} \vskip -5mm \indent}
\def\Sec{\@Startsection{section}{1}{\z@}
                                   {-3.5ex \@plus -1ex \@minus -.2ex}%
                                   {2.3ex \@plus.2ex}%
                                   {\normalfont\large\bfseries\boldmath}} % czj
\def\@Startsection#1#2#3#4#5#6{%
  \if@noskipsec \leavevmode \fi
  \par
  \@tempskipa #4\relax
  \@afterindenttrue
  \ifdim \@tempskipa <\z@
    \@tempskipa -\@tempskipa \@afterindentfalse
  \fi
  \if@nobreak
    \everypar{}%
  \else
    \addpenalty\@secpenalty\addvspace\@tempskipa
  \fi
  \@ifstar
    {\@ssect{#3}{#4}{#5}{#6}}%
    {\@dblarg{\@Sect{#1}{#2}{#3}{#4}{#5}{#6}}}}
\def\@Sect#1#2#3#4#5#6[#7]#8{%
  \ifnum #2>\c@secnumdepth
    \let\@svsec\@empty
  \else
    \refstepcounter{#1}%
    \protected@edef\@svsec{\@seccntformat{#1}\relax}%
  \fi
  \@tempskipa #5\relax
  \ifdim \@tempskipa>\z@
    \begingroup
      #6{%
          \begin{center} % czj
          \@hangfrom{\hskip #3\relax\S$\,$\@svsec\hskip -4mm . \hskip 1mm}%
          \interlinepenalty \@M #8\@@par 
          \end{center}}%
    \endgroup
    \csname #1mark\endcsname{#7}%
    \addcontentsline{toc}{#1}{%
      \ifnum #2>\c@secnumdepth \else
        \protect\numberline{\csname the#1\endcsname}%
      \fi
      #7}%
  \else
    \def\@svsechd{%
      #6{\hskip #3\relax
      \@svsec #8}%
      \csname #1mark\endcsname{#7}%
      \addcontentsline{toc}{#1}{%
        \ifnum #2>\c@secnumdepth \else
          \protect\numberline{\csname the#1\endcsname}%
        \fi
        #7}}%
  \fi
  \@xsect{#5}}
\def\list#1#2{\ifnum \@listdepth >5\relax \@toodeep \else \global
\advance \@listdepth\@ne \fi \rightmargin \z@ \listparindent\z@
\itemindent\z@ \csname @list\romannumeral\the\@listdepth\endcsname
\def\@itemlabel{#1}\let\makelabel\@mklab \@nmbrlistfalse #2\relax
\@trivlist \parskip 0pt \parindent\listparindent \advance \linewidth
-\rightmargin \advance\linewidth -\leftmargin \advance\@totalleftmargin
\leftmargin \parshape \@ne \@totalleftmargin \linewidth \ignorespaces}
\catcode`@=12
\theoremstyle{plain}
\newtheorem{thm}{\indent Theorem}[section]
\newtheorem{lem}{\indent Lemma}[section]
\newtheorem{prop}{\indent Proposition}[section]

\theoremstyle{remark}

\def\thebibliography#1{\section*{\centerline{\bf \large References}}
\list{[\arabic{enumi}]} {\settowidth \labelwidth{[#1]} \leftmargin
\labelwidth \advance \leftmargin \labelsep \usecounter{enumi}}
\def\newblock{\hskip .11em plus .33em minus .07em} \footnotesize \sloppy \clubpenalty
4000 \widowpenalty 4000 \sfcode`\.=1000 \relax}

\newcommand{\bthmm}{\begin{thm}}
\newcommand{\ethmm}{\end{thm}}
\newcommand{\blemm}{\begin{lem}}
\newcommand{\elemm}{\end{lem}}
\newcommand{\bppm}{\begin{prop}}
\newcommand{\eppm}{\end{prop}}
\newcommand{\bprfm}{\begin{Proof}}
\newcommand{\eprfm}{\end{Proof}}
\newcommand{\bcorm}{\begin{cor}}
\newcommand{\ecorm}{\end{cor}}
\newtheorem{theorem}{Theorem}[section]{\bf}{\it}
\newtheorem{proposition}{Proposition}[section]{\bf}{\it}
\newtheorem{lemma}{Lemma}[section]{\bf}{\it}
{\it}{\rm}
\newtheorem{remark}{Remark}[section]{\it}{\rm}
\newtheorem*{Proof}{Proof}{\it}{\rm}
\newcommand{\bthm}{\begin{theorem}}
\newcommand{\ethm}{\end{theorem}}
\newcommand{\blem}{\begin{lemma}}
\newcommand{\elem}{\end{lemma}}
\newcommand{\bpp}{\begin{proposition}}
\newcommand{\epp}{\end{proposition}}
\newcommand{\bprf}{\begin{proof}}
\newcommand{\eprf}{\end{proof}}
\def\nd{\noindent}
\def\hs#1{\hskip 2pt #1\hskip 2pt}

\newcommand{\fa}{\forall}

\newcommand{\q}{\quad}
\newcommand{\qq}{\qquad}

\newcommand{\fr}{\frac}
\newcommand{\sq}{\sqrt}

\newcommand{\s}{\ }

\newcommand{\tm}{\times}

\newcommand{\vep}{\varepsilon}

\newcommand{\gm}{\gamma}

\newcommand{\al}{\alpha}
\newcommand{\bt}{\beta}
\newcommand{\tta}{\theta}

\newcommand{\beq}{\begin{eqnarray}}
\newcommand{\eeq}{\end{eqnarray}}
\newcommand{\beqs}{\begin{eqnarray*}}
\newcommand{\eeqs}{\end{eqnarray*}}
\newcommand{\bc}{\begin{center}}
\newcommand{\ec}{\end{center}}
\newcommand{\cl}{\centerline}
\newcommand{\bcase}{\begin{cases}}
\newcommand{\ecase}{\end{cases}}
\newcommand{\bmat}{\begin{matrix}}
\newcommand{\emat}{\end{matrix}}
\newcommand{\bbm}{\begin{bmatrix}}
\newcommand{\ebm}{\end{bmatrix}}
\newcommand{\bpm}{\begin{pmatrix}}
\newcommand{\epm}{\end{pmatrix}}
\newcommand{\bvm}{\begin{vmatrix}}
\newcommand{\evm}{\end{vmatrix}}
\newcommand{\bal}{\begin{align*}}
\newcommand{\eal}{\end{align*}}
\newcommand{\bale}{\begin{aligned}}
\newcommand{\eale}{\end{aligned}}
\newcommand{\bequ}{\begin{equation}}
\newcommand{\eequ}{\end{equation}}

\newcommand{\mf}{\mathfrak} %Germany
\newcommand{\mbbr}{\mathbb{R}}
\newcommand{\mbbc}{\mathbb{C}}
\newcommand{\mfhp}{\mathfrak{h}_{\mathfrak{p}}}
\newcommand{\mfa}{\mathfrak{a}}
\newcommand{\mfb}{\mathfrak{b}}
\newcommand{\mfc}{\mathfrak{c}}
\newcommand{\mfd}{\mathfrak{d}}
\newcommand{\mfe}{\mathfrak{e}}
\newcommand{\mff}{\mathfrak{f}}
\newcommand{\mfg}{\mathfrak{g}}
\newcommand{\mfh}{\mathfrak{h}}

\newcommand{\mfk}{\mathfrak{k}}
\newcommand{\mfl}{\mathfrak{l}}

\newcommand{\mfp}{\mathfrak{p}}

\newcommand{\mfu}{\mathfrak{u}}

\newcommand{\mfz}{\mathfrak{z}}
\newcommand{\bpic}{\begin{picture}}
\newcommand{\epic}{\end{picture}}

\renewcommand{\thefootnote}{\fnsymbol{footnote}}

\title{\Large \bf \boldmath CURVATURE ESTIMATES FOR IRREDUCIBLE SYMMETRIC SPACES}

\author{\sc \large Liu Xusheng$^\ast$}

\date{}

\begin{document}

\maketitle
\thispagestyle{first}
\renewcommand{\thefootnote}{\fnsymbol{footnote}}

\renewcommand{\thefootnote}{\fnsymbol{footnote}}
\footnotetext{\hspace*{-5mm}
 \begin{tabular}{@{}r@{}p{13.4cm}@{}}
& Manuscript received  \\
$^\ast$ & Department of Mathematics,Fudan University, Shanghai, China.\\
&{\bf E-mail:xshliu@fudan.edu.cn } \\
\end{tabular}}

\renewcommand{\thefootnote}{\arabic{footnote}}

\begin{abstract}
By making use of the  classification of real simple Lie algebra,
we get the maximum of the squared length  of  restricted roots
case by case,   thus we get the upper bounds of sectional curvature for irreducible
Riemannian symmetric spaces of compact type.
As an application,   we verify Sampson's  conjecture
in  all cases for irreducible  Riemannian symmetric spaces of noncompact type.

\vskip 4.5mm

\nd \begin{tabular}{@{}l@{ }p{10.1cm}} {\bf Keywords } & Symmetric
space, Semi-simple Lie algebra, Harmonic map.
\end{tabular}\\
\nd {\bf 2000 MR Subject Classification }  22E60, 53C35, 53C43 \\
\nd {\bf Chinese Library Classification} O152.5, O186.12 \\
\nd {\bf Document Code} A \\

\end{abstract}
\baselineskip 14pt
\setlength{\parindent}{1.5em}

\setcounter{section}{0}
\Section{Introduction and statement of results}
\vskip -3mm

There are many Liouville type theorems for harmonic maps. It was
conjectured by J. H. Sampson [9] that any harmonic map with finite
energy from a complete simply connected Riemannian manifold  with
non positive sectional curvature whose dimension exceeds 2 must be
constant. This is valid  for space forms,   but unsolved in
general case. For Cartan-Hadamard manifolds,   Y. L. Xin [11]
proved a general vanishing theorem as follows.  \par \bthm
\label{thma} Let M be an m-dimensional Cartan-Hadamard manifold
with the sectional curvature  $-a^2 \le K\le 0$ and the Ricci
curvature bounded from above by $-b^2$.  Let f be a harmonic map
from M into any Riemannian manifold with the moderate divergent
energy. If $b\ge 2a$,   then f has to be constant. \ethm
\par
Moreover, by using $L^2$-cohomology technics, from  Theorem 2.2 of [5],
the above condition  can be relaxed to $b\ge \sq 2 a$.

For the  irreducible symmetric space,   the Ricci curvature is
constant. In this article we get the bounds of sectional curvature
by using the method of Lie  algebra. As an application,  we verify
Sampson's  conjecture in  most cases for Riemannian symmetric
spaces of noncompact type by using the above Theorem 1.1.
\par

\bthm \label{thmb}  Let M be one of the irreducible symmetric spaces of
noncompact type in the following cases(see [4,  p518,  Table V])
\beqs
&& SL(n,  \mbbr)/SO(n),  n\ge 8; \q
SU^*(2n)/Sp(n),  n\ge 2; \\
&& SU(p, q)/S(U_p\tm U_q),  p+q\ge 8; \\
&&SO_o(p, q)/SO(p)\tm SO(q), p+q\ge 6, for \s r=1, \\
&& \hspace{2cm} p+q\ge 10 \s for\s  r>1,  where\s  r=min(p, q); \\
&& SO^*(2n)/U(n),  n\ge 5; \q Sp(n, \mbbr)/U(n),  n\ge 7; \\
&&Sp(p, q)/Sp(p)\tm Sp(q),  p+q\ge 3; \\
&&  EI, EII, EIII, EIV,  EV, EVI, EVII, EVIII,  EIX,  FI,FII \s and \s G.
\eeqs
Let f be a harmonic map from M into
any Riemannian manifold with the moderate divergent energy,
then f has to be constant.
\ethm
\par

In \S 2,  we collect basic  facts about real simple Lie algebras
and symmetric spaces. In \S 3,  we give the upper bounds of
sectional curvature  for Riemannian symmetric spaces of compact
type by choosing the Cartan subalgebra with maximal vector part.
In \S 4,  we give the method of obtaining  Cartan subalgebra with
maximal vector part from a Cartan subalgebra with maximal torus
part. In \S 5,  we carry out the detailed calculation on estimates
of sectional curvature by using the root space decompositions of
complex   simple Lie algebras. The final results will be shown in
Table 5.1. This is our main result.  It is interesting in its own
right, although some of the classical cases were given by Y.C.Wong
[10].

From Table 5.1 and Theorem \ref{thma}, Theorem \ref{thmb} follows
immediately.

\Section{Some basic facts for real semisimple Lie algebra \\
and irreducible symmetric space}
\vskip -3mm

Let $M=U/K$ be a compact irreducible Riemannian symmetric space of type I,
its noncompact dual is $G/K$, which  is of type III.
There are decompositions for corresponding Lie algebras
\beq
\mf{u}=\mf k+i\mfp, \q \mfg_0=\mfk + \mfp, \q [\mfk, \mfk]\subset \mfk,\q  [\mfk, \mfp]\subset \mfp, \q [\mfp, \mfp]=\mfk
\eeq
Let $\mfg=\mfu^{\mbbc}=\mfg_0^\mbbc$ be its complexified simple Lie algebra,
the rank of $\mf g$ (resp.  M) is  the dimension of
the maximal abelian complex (resp.  real) subspace.
The Killing form of $\mfg$ is
\beq (X, Y)=B(X, Y)=tr ad X ad Y. \eeq
Let
\beq \label{innp}
<X, Y>=\epsilon (X, Y)
\eeq
where $\epsilon=1$ for noncompact type;  -1 for compact type.
The restriction on $i\mfp$(resp.  $\mfp$)
of $<, >$ gives the Riemannian  inner product at $o=\pi(e)$ for the corresponding
symmetric space, where $e$ is the identity of U.
The curvature tensor is $R(X, Y)=-ad[X, Y]$,  from the invariance
of Killing form  we get
\beq
R(X, Y, Z, W)=\epsilon (-[[X, Y], W], Z)=\epsilon ([X, Y], [Z, W]).
\eeq
Since the Killing form is positive definite on $\mfp$
and negative definite on $i\mfp$,
the curvature is  nonnegative for symmetric space of compact type,
non positive for  symmetric space of noncompact type.
Moreover,  by choosing a  adapted base  we have ([8, p180])
\beq
Ric(X, Y)=-\fr 12 (X, Y).
\eeq
So the Ricci curvature is
$\fr 12$ for compact type and $-\fr 12$ for noncompact type.
Owing to the dual relation,
we only need to calculate the upper bounds of sectional curvature
for the irreducible Riemannian symmetric spaces of compact type,
and the lower bounds are 0 except for the case of rank 1.
\par

The irreducible  Riemannian symmetric spaces of type II and  IV  are dual to each other.
Let M be of type II,  it can be identified as a connected  compact simple  Lie group U with the bi-invariant
metric,  see [4,  p516].
We give the Riemannian inner product at e
\beq \label{cpt00}
<X, Y>=-(X, Y).
\eeq
We have the Cartan decomposition of $\mfg=\mfu^\mbbc$
and the decomposition of  root spaces.
For $X, Y\in \mfu$,  we have ([8,  p185])
\beq
\nabla_XY=\fr 12[X, Y].
\eeq
So
\beq
R(X, Y)=-\fr 14ad[X, Y].
\eeq
In next section,  we get that the Ricci curvature is $\fr 14$.

Recall that the pair $(\mfl,  \theta)$ is called effective orthogonal
symmetric Lie algebra if (i)$\mfl$ is a real Lie algebra;
(ii)$\tta$ is an involutive automorphism of $\mfl$, i.e., $\tta \ne id$ but $\tta^2=id$;
(iii) the set of fixed points of $\tta,  \mfk, $  is a compactly
imbedded subalgebra of $\mfl$;
(iv)$\mfk\cap \mfz(\mfl)=0$,  where $\mfz(\mfl)$ is the center of $\mfl$.
The pair $(G, K)$ is said to be associated with $(\mfl, \tta)$  if
G is a connected Lie group with Lie algebra $\mfl$,  and K is a Lie
subgroup of G with Lie algebra $\mfk$.
If K is closed in G,  $G/K$ endowed with the G-invariant metric
becomes a Riemannian symmetric space.
The orthogonal symmetric Lie algebra $(\mfl,  \theta)$ is called
irreducible if,  let $\mfl=\mfk+\mfp$ be the $\pm 1$ characteristic
subspace decomposition of $\tta$,  $\mfl$ be real semisimple Lie algebra,
$\mfk$ contains no nonzero idea of $\mfl$,
and the adjoint algebra $ad_\mfl \mfk$ acts irreducibly on $\mfp$.
If $\mfl$ is a compact simple real Lie algebra,  $(\mfl, \tta)$ is of type I,
the associated Riemannian symmetric space is of type I;
if $\mfl$ is a noncompact real simple Lie algebra,
its complexification $\mfl^\mbbc$ is a complex simple Lie algebra,  then
$(\mfl, \tta)$ is of type III,
the associated Riemannian symmetric space is of type III.
Furthermore,  every  irreducible symmetric Lie algebra of type I corresponds
1-1 to one of type III,  they are dual to each other. \par

A real simple Lie algebra is of first kind if its complexification
is a complex simple Lie algebra;  otherwise of second kind.
Let $\mfg$ be a complex semisimple Lie algebra,  there exists a compact real form $\mfu$,
and for any two compact real forms,  there exists an inner automorphism
of $\mfg$ which transforms  one to the other,
which means there is only one compact real form up to equivalence.
\par
A maximal abelian subalgebra of a complex semisimple Lie algebra
is called a Cartan subalgebra of $\mfg$,
any two Cartan subalgebras of $\mfg$ are conjugate
under the group of inner automorphism.
A maximal real abelian subalgebra of a real semisimple Lie algebra $\mfg_0$
is called a Cartan subalgebra of $\mfg_0$.
 The Cartan subalgebra $\mfh_0$ of $\mfg_0$ is a real form of a Cartan  subalgebra $\mfh$
 of $\mfg=\mfg_0^\mbbc$,   namely $\mfh_0=\mfh\cap \mfg_0$.  \par

The real simple Lie algebra of first kind can be classified as follows.
Given a complex simple Lie algebra $\mfg$,  let $\mfu$ be its compact real form,
let $\theta$ be an involution which leaves $\mfu$ invariant,
the Cartan decomposition is
$\mfu=\mfk+i\mfp$,  where $\mfk, i\mfp$ are $\pm 1$ characteristic
subspace of $\theta$,
$[ \mfk,  \mfp]\subset \mfp, [ \mfp,  \mfp]\subset \mfk, \mfk$
is also a subalgebra of $\mfu$,  called  characteristic subalgebra with respect to $\tta$.
Let $\mfg_0=\mfk+\mfp$,
if $\tta\ne 1$, then $\mfg_0$ is a noncompact real simple Lie algebra.
Two involutions determine the same real simple Lie algebra (up to an automorphism)
if and only if they are conjugate under an inner automorphism.
So every noncompact real simple Lie algebra corresponds 1-1 to a conjugacy class
of involution of $\mfg$,
and every  conjugacy class corresponds 1-1  to a Riemannian symmetric space of compact type (or of noncompact type).
Here $\tta$ is the involution of complex simple Lie algebra $\mfg$,
meanwhile it is  the involution of real simple Lie algebra $\mfg_0$ and $\mfu$.
For details referring to [3], [7]. \par

Let $\mfh$ be a Cartan subalgebra of $\mfg$,  with respect to $\mfh$ there is the
root space decomposition $\mfg=\mfh+\sum_{\al \in \Delta} \mfg_{\al}$.
Denote the conjugation of $\mfg$ with respect to $\mfu$  (resp.  $\mfg_0$) by $\tau$  (resp. $\sigma$).
Define
\beq \label{inp}
<X, Y>=-(X,  \tau Y), \q |X|=\sq{<X, X>}
\eeq
This gives  a Hermitian inner product on $\mfg$.
Any root $\al$ can be imbedded into $\mfh$ by
\beq \label{rooti}
\al (h)=(\al,  h),  h\in \mfh.
\eeq
We call the real subspace $\mfh_0=\sum_{\al\in \Delta} \mbbr \alpha$  generated by all roots
the real part of $\mfh$.
Since the root is purely imaginary with respect to $\mfu,  \al\in i\mfu$,
there exists the decomposition
\beq
\mfh_0=i\mfh_{\mfk}+\mfh_{\mfp},  \mfh_{\mfk}=i\mfh_0\cap\mfk, \mfh_{\mfp}=\mfh_0\cap \mfp.
\eeq
\par
Choose the Weyl basis
\beq
&e_{-\al}=\tau (e_\al),  (e_\al,  e_{-\al})=(e_\al,  \tau(e_\al))=-1\\
&[e_\al,  e_\bt]=N_{\al\bt}e_{\al\bt},  N_{\al\bt}=-N_{\bt\al}, N_{-\al -\bt}=N_{\al\bt}\\
&[e_\al, e_\al]=-\al
\eeq
where $N_{\al\bt}^2=\fr 12 q_{\bt\al}(1+p_{\bt\al})(\al, \al)\in \mbbr,
 \bt -p_{\bt\al}\al,  \cdots,  \bt+q_{\bt\al}\al$ is the $\al$-chain of $\bt$. \par
$\mfu$ is generated over $\mbbr$ by
\beq \label{uv} \s\s \{i\alpha,  u_\alpha,  v_\alpha,  \alpha \in \Delta^+ \},
 u_\alpha=\fr 1{\sq 2}(e_\alpha+e_{-\alpha}),
 v_\alpha=\fr {i}{\sq 2}(e_\alpha-e_{-\alpha}).
\eeq
Any element of $\mfu$ can be expressed as
 $$X=ih+\sum_{\al\in \Delta^+}a_{\al}e_{\al}+\bar{a}_{\al}e_{-\al} \qq h\in \mfh_0,  a_{\al}\in \mbbc. $$
Moreover
\beq
 [ih,  u_\alpha]=(\al, h)v_\alpha, \qq
 [ih,  v_\alpha]=-(\al, h)u_\al, \qq h\in \mfh_0.
\eeq
 The Cartan subalgebra of $\mfk$ is called a torus Cartan subalgebra of $\mfg_0$.
A Cartan subalgebra of $\mfg_0$ is called a Cartan subalgebra with maximal torus part
if it contains a torus Cartan subalgebra.
A maximal abelian subspace of $\mfp$ is called a  reduced Cartan subalgebra of $(\mfu,  \tta)$,
a Cartan subalgebra of $\mfg_0$ is called a  Cartan subalgebra with maximal vector part
if it contains a reduced Cartan subalgebra.  \par

In general,  we can choose Cartan subalgebra in different ways.
The real simple Lie algebra of first kind can be given as follows.
Let $\mfh_T$ be a Cartan subalgebra of $\mfk$,
it can be extended to a  Cartan subalgebra of $\mfg$,  say  $\mfh$,
then $\mfh$ is a Cartan subalgebra of $\mfg$ with maximal torus part.
Let $\Delta$ be  the root system of $\mfg$ with respect to $\mfh$,
we  can choose the compatible order,
let
\beq
\Pi=\{\al_1,  \al_2, \cdots, \al_l\}
\eeq
be a simple root system,
then $\tta$ has Gantmacher standard form ([2], [6]):
\beq \label{gan1}
\theta=\theta_0exp(2\pi i ad h_0), h_0\in i\mfh_T
\eeq
where $\theta_0$ is a special rotation with respect to $\Pi$,
$\theta_0$ preserves $\Pi$ and be an involution which called a regular involution,
$\tta_0(h_0)=h_0$. \par
If $\tta_0=1$,  denote the highest root of $\Delta$ by
$$\delta=\sum_{i=1}^l m_i(\delta)\al_i$$
we can choose $i,  m_i(\delta)=1$ or 2 such that
\beq  \label{gan2}
(h_0,  \al_i)=\fr 12,  \qq  (h_0,  \al_j)=0,  j\ne i.
\eeq
\par
If $\tta_0\ne 1$,  let the characteristic subalgebra of $\tta_0$ be $\mfk_0=\{X\in \mfu | \tta_0(X)=X\}$,
the Cartan subalgebra $\mfh_0'$ of $\mfk_0$  is the intersection of
a Cartan subalgebra  $\mfh_0$ of $\mfg_0$ with $\mfk_0$,
$\mfh_0'=\mfh_0\cap \mfk_0$,  let $\{\al_i'\}$ be the simple root system of $\mfh_0'$, let
$$\lambda'=\sum_i m_i(\lambda')\al_i'$$ be the highest root of $\mfk_0'$,
then we can choose i,  $m_i(\lambda')=1 $ or 2  so that
\beq  \label{gan3}
(h_0,  \al_i')=\fr 12,  \qq  (h_0,  \al_j')=0,  j\ne i.
\eeq
\par
For these two cases we denote $h_0$ by $h_i$.
So for a given  complex simple Lie algebra,
by  finding  all $\tta_0,  h_i$,
we can exhaust all real simple Lie algebras of first kind.
$\theta_0\ne 1$ only when $\mfg=\mfa_l,  \mfd_l, \mfe_6$.
For details referring to [7].

\Section {The Cartan subalgebra with maximal vector part \\
and the upper bounds of sectional curvature}
\vskip -3mm

Let $M=U/K$ be an irreducible Riemannian symmetric space of compact type,
let $\mfh_{\mfp}$ be  a maximal abelian subspace of $\mfp$,
it can be extended to a Cartan algebra $\tilde\mfh$ of $\mfg$
by choosing  $\mfh_{\mfk}\subset \mfk$,
\beq \tilde \mfh_0=i\mfh_{\mfk}+\mfh_{\mfp},  \tilde\mfh=\mfh_0^C.  \eeq
\par
The root system of $\mfg$ with respect to $\tilde\mfh$ is $\Delta$,
the simple root system is
\beq
\Pi=\{\al_1,  \al_2, \cdots, \al_r,  \cdots,  \al_l\}
\eeq
where $\al_1,  \al_2, \cdots, \al_r\in \mfh_\mfp$.
\par
The tangent space of M is generated by
\beq \{i\al_1,  \cdots,  i\al_r,  u_\alpha-\theta(u_\al),  v_\alpha-\theta(v_\al),  \al\in \Delta^+ \} \eeq
where $u_\al,  v_\al$ is defined by (\ref{uv}). \par

The projection on $\mfh_{\mfp}$ of a root $\al$ is $\al'=\fr 12(\al-\theta (\al))$,
if $\al'\ne 0$,  $\al'$ is called a restricted root.
In the same way by $(h,  \al')=\al'(h),  h\in \mfh_{\mfp}$,  any
restricted root can be identified with an element in $\mfh_{\mfp}$.
It's well known that every restricted root is  the restriction
of a root.  All restricted roots span  $\mfh_{\mfp}$ over $\mbbr$. We also have
\beq
[ h, e_{\al} ]=(\al', h)e_{\al}, h\in \mfh_{\mfp}.
\eeq
Denote the set of restricted roots by $\Sigma'=\{\al',  \al'\ne 0\}$,
the set of noncompact roots by $\Sigma=\{\al\in \Delta,  e_\al\in \mfp^\mbbc  \}$. \par

We remark that there are similar  projections  by choosing different Cartan subalgebra.
In the sequel  we choose the Cartan subalgebra with maximal torus part.
By abuse of notation,  we denote
the projection on torus part or vector part by the same symbol $\al'$,
for torus part $\al'=\fr {\al+ \tta(\al)}2$,
but for vector part ,  $\al'=\fr {\al- \tta(\al)}2$.
\par
From the general theory of symmetric space we know that
any tangent vector is conjugate to  an element  $ih$ of $i\mfh_{\mfp}$
 by the action of isotropy group([4, p246] Theorem 6.2).
So it is only to estimate the sectional curvature of the form
$$K(ih,  X)\qq  h\in \mfhp,  X\perp \mfhp, X\in \mfp$$
Then
\beq
X=\sum_{\al\in \Sigma^+} a_\al e_\al+\bar a _{\al} e_{-\al}.
\eeq
\beq
[ ih,  X ]=\sum i(\al,  h)(a_\al e_\al-\bar a _{\al} e_{-\al})
\eeq
\beq
R(ih, X, ih, X)& = &|[ ih,  X ] |^2=\sum 2(h, \al )^2|a_\al|^2 \\
&=&\sum 2(h, \al ')^2|a_\al|^2\le |h|^2 \sum 2(\al', \al')^2|a_\al|^2.
\eeq
If $|h|=1,  |X|=1$ then $\sum |a_\al|^2=\fr 12$ ,  let d be the maximum of
squared length of restricted roots, then
$$K(ih, X)\le d. $$
On the other hand,  if $|\al_0'|^2=d$, \s set $h=\al_0', X=e_{\al_0}+e_{-\al_0}$
it can be verified that
\beq
K(ih, X)=(\al_o', \al_0)= d.
\eeq
So we get the following proposition ([4, p334] Theorem 11.1):\par
\bpp
The maximum of sectional curvature for a irreducible Riemannian
symmetric space of compact type is the squared length of the highest restricted root.
\epp
\par
Let $\delta$ be highest noncompact root,  then $\delta'$ is  highest  restricted root.
So in order to get the maximum of sectional curvature it is sufficient to
get the maximum of the squared length of the projection of noncompact roots
onto the reduced Cartan subalgebra. \par

Now we proceed to the  irreducible  Riemannian symmetric space of type II.
Let $M=U$ be a compact simple Lie group with Lie algebra  $\mfu$,
the Riemannian metric
is given in (\ref{cpt00}), let $\mfh$ be a Cartan subalgebra
of $\mfg=\mfu^\mbbc$,  the root system of $\mfg$ with respect to $\mfh$
is $\Delta$. For $h\in \mfh_0$,
\beqs
(h, h)&=& tr \   ad h \cdot ad h\\
&=& \sum_{\al\in \Delta}(\al, h)^2\\
&=& 2\sum_{\al>0}(\al, h)^2.\\
Ric(ih, ih)&=&\sum_{\al>0}R(ih,  u_\al,  ih,  u_\al)+R(ih,  v_\al,  ih,  v_\al)\\
&=&-\fr 14\sum_{\al>0}([ih,  u_\al],  [ih,  u_\al])+([ih,  v_\al],  [ih,  v_\al])\\
&=& -\fr 14\sum_{\al>0}(\al, h)^2(v_\al,  v_\al)+(\al, h)^2(-u_\al,  -u_\al)\\
&=& \fr 12\sum_{\al>0}(\al, h)^2\\
&=& -\fr 14(ih, ih)\\
&=&\fr 14<ih, ih>.
\eeqs
We see that the Ricci  curvature is $\fr 14$.
By the same method, we get that  the upper bound of curvature for U is exactly
the maximum of the squared length of roots,
ie.,  the squared length of the highest root.

 \par
In fact, let $U^*=\{(u, u)|u\in U\},  \mu (u, v)=(v, u)$,
let $U\tm U/U^*$
be the  Riemannian symmetric pair of compact type associated with
the orthogonal symmetric Lie algebra
$(\mfu\tm \mfu,  d\mu)$.
This metric    coincides with
the previous one in (\ref{cpt00})  up to a constant.

\Section {The transformation from the Cartan subalgebra with maximal torus part to
  the Cartan subalgebra with maximal vector part}
\vskip -3mm

Let $\mfg$ be complex simple Lie algebra with the Cartan decomposition
\beq
\mfu=\mfk +i\mfp,  \mfg_0=\mfk+\mfp.
\eeq
Let $\mfh_0=i\mfh_T+\mfh_{V}, \mfh_T\subset \mfk,  \mfh_V\subset \mfp,  \mfh=\mfh_0^C$
be a Cartan subalgebra with  maximal torus part,  the  root system is  $\Delta$,
$\theta=\theta_0exp(2\pi i ad h_0) $ as in (\ref{gan1}) - (\ref{gan3}). \\
The decomposition of root system is
\beq
&\Delta=\Delta_0\cup \Delta_n\cup\Delta_c\\
&\Delta_0=\{\al\in \Delta|\theta(\al) \ne \al  \}\\
&\Delta_n=\{\al\in \Delta|\theta(\al)=\al,  \theta(e_\al)=-e_\al \}\\
&\Delta_c=\{\al\in \Delta|\theta(\al)=\al,  \theta(e_\al)=e_\al \}.
\eeq
The centralizer of $\mfh_V^C$ in $\mfg$ is $\mfz(\mfh_V)=\{X\in g| [X, h]=0, h\in \mfh_V\}$,
we have
\beq
\mfz(\mfh_V)=\mfh_V^C+\sum_{\al\in \Delta_n}\mfg_{\al}.
\eeq
We call a root $\al$ compact if $\mfg_\al\subset \mfk^\mbbc$,  otherwise noncompact.  \par
We now give  the method of Harish-Chandra for constructing
Cartan subalgebra  with maximal vector part
from a Cartan subalgebra with maximal torus part.
It is so-called Cayley's transformation.
For details referring to [12],   see also in  [4, p385-387, p530-534].
\par
Two roots $\al,  \bt$ are called strongly orthogonal if $\al\pm\bt\notin \Delta$. \par
Let $\gm_1$ be the highest root in $\Delta_n^+$,  set
\beq
\Delta_n^+(\gm_1)=\{\al\in \Delta^+_n|\al\ne \gm_1,  \al\pm\gm_1\notin \Delta\}.
\eeq
Let $\gm_2$ is the highest root in $\Delta_n^+(\gm_1)$,  set
\beq
\Delta_n^+(\gm_1, \gm_2)=\{\al\in \Delta^+_n(\gm_1)|\al\ne \gm_2,  \al\pm\gm_2\notin \Delta\}.
\eeq
In this way we can get the maximal system of positive roots which are strongly orthogonal to
each other,  say $\{\gm_1, \gm_2, \cdots,  \gm_s\}$. \par
Let $\mfa_1$ be the real subspace generated by $\{e_{\gm_i}-e_{-\gm_i},  i=1, 2, \cdots, s\}$,
then
\beq
\mfa=\mfa_1+\mfh_V
\eeq
is a reduced Cartan subalgebra of $\mfp$. \par
Let
\beq
\tilde\mfh_T=\{X\in \mfh_T|(X,  \gm_i)=0, i=1, 2, \cdots, s \}
\eeq
then
\beq
\tilde\mfh_0=\tilde\mfh_T+\mfa
\eeq
is a  Cartan subalgebra with  maximal vector part of $\mfg_0$. \par
Let
\beq
&X_i=\fr {\pi}{2\sq 2 |\gm_i|}(e_{\gm_i}+e_{-\gm_i})\\
&\rho_i=exp(adX_i)\in Int(\mfu)\\
&\rho=\rho_1\rho_2\cdots\rho_s
\eeq
then
\beq
&\rho(\mfh_0)=\tilde\mfh_0\\
&\rho(\gm_i)=-\fr{|\gm_i|}{\sq 2}(e_{\gm_i}-e_{-\gm_i})\\
&\rho|\tilde\mfh_T+\mfh_V=id.
\eeq
\par
Being an  automorphism,   $\rho$ exchanges two Cartan subalgebra,
$\rho$ take  roots to roots,  eigenvectors to eigenvectors.
Let $\al_0$ be a noncompact root which takes maximal root length,
Since $\rho$ is a isometry of $\mfu$,
 The projections on  two Cartan subalgebra $\mfh_0, \tilde\mfh_0$
 of $\al_0, \rho(\al_0)$  have the same  length.
 We need only to calculate the maximal length
 of the projection on $\{\gm_1, \gm_2, \cdots, \gm_s\}\cup \mfh_V$
 of non compact roots with respect to $\mfh$.
 This is the maximal  length of restricted roots with respect
 to $\tilde\mfh_0$,  the  Cartan subalgebra with maximal vector part.
\par
We can choose compatible order so that $\rho$ take the highest noncompact
root in  $\mfh_0$ to the highest noncompact root in $\tilde\mfh_0$.
 Let $\al_0$ be the highest noncompact root with respect to $\mfh_0$
 ,  the projection  on $\{\gm_1, \gm_2, \cdots, \gm_s\}\cup \mfh_V$ is $\al_0'$,
then $|\al_0'|$ attains the maximal  length of restricted roots. \par

Since real simple Lie algebras of the first kind are completely classified
according to Cartan subalgebra with maximal torus part,
the irreducible symmetric spaces of type I and type III
can be given by using the Gantmacher standard form of complex
simple Lie algebras.
Moreover by estimating  the maximum of squared length of restricted roots
we get the curvature bounds for symmetric spaces.

\Section {estimates of bounds of sectional curvature}
\vskip -3mm

Let $\mfg$ be a complex simple Lie algebra,
let $\mfu=\mfk +i\mfp$ be the Cartan decomposition of $U/K$ which is a irreducible symmetric
space of compact type,
as in (\ref{inp}),  the Riemannian inner product  at $\pi(e)$ is
\beq
<X, Y>=-(X, Y).
\eeq
Every root is embedded in $\mfg$ through Killing form as in (\ref{rooti}).
Then $\fa \al\in \Delta,  \al\in i\mfu,  (\al, \al)>0$.
For every simple Lie algebra,  there are at most two different lengths of roots.
Let $\mfh$ be a Cartan subalgebra,  for $H_1,  H_2\in \mfh$,
\beq \q
(H_1, H_2)=tr (ad H_1 ad H_2)=\sum_{\bt \in \Delta} \bt(H_1)\bt(H_2)
=2\sum_{\bt \in \Delta^+} \bt(H_1)\bt(H_2).
\eeq
We have
\beq\q
(\al, \al)&=&2\sum_{\bt \in \Delta^+} (\bt, \al)(\bt, \al)=
\fr 12 \sum_{\bt \in \Delta^+}(\fr {2(\bt, \al)}{(\al, \al)})^2(\al, \al)^2\\
&=&\fr 12 \sum_{\bt \in \Delta^+} a_{\bt\al}^2(\al, \al)^2
\eeq
where
$a_{\bt\al}=\fr {2(\bt, \al)}{(\al, \al)}=p_{\bt\al}-q_{\bt\al}$ be Cartan integer,
so
\beq \label{ci}
(\al, \al)=\fr 2{\sum_{\bt \in \Delta^+} a_{\bt\al}^2}.
\eeq
By this normalization we can calculate the maximum of squared length
of roots, as listed  in [1].

We adopt the  Dynkin diagrams of complex simple Lie algebras  as in [4, p476]. \par

For the Gantamacher standard form as in (\ref{gan1}),
$\tta_0\ne 1$ only when $\mfg=\mfa_l,  \mfd_l,  \mfe_6$.
We list these special rotations as follows
(for Dynkin diagrams referring to the following paragraphs).
\par
For $\mfa_l,  \tta_0(\al_1)=\al_{l},   \tta_0(\al_2)=\al_{l-1}, \cdots,   \tta_0(\al_l)=\al_{1}$\par
For $\mfd_l,  \tta_0(\al_{l-1})=\al_{l},   \tta_0(\al_l)=\al_{l-1}, \tta_0(\al_j)=\al_{j},  j\le l-2$\par
For $\mfe_6,  \tta_0(\al_{1})=\al_{6},   \tta_0(\al_2)=\al_{2}, \tta_0(\al_3)=\al_{5},
\tta_0(\al_4)=\al_{4}$,  \\
\hspace*{1.5cm}$\tta_0(\al_5)=\al_{3}, \tta_0(\al_6)=\al_{1}.$
\par
We can reduce the calculation to the  following cases:\par
(1)If $\theta=exp(2\pi i ad h_i)$ is an inner automorphism,  from the classification of
real simple Lie algebra we have $\mfg_0=\mfk+\mfp,  \mfh_0=i\mfh_{T}$ which means
$\mfk$ is a maximal compact subalgebra  of $\mfg_0$,
its Cartan subalgebra is also a Cartan subalgebra of $\mfg_0$,
$\mfa$ consists of r strongly orthogonal roots in $\Delta_n^+$.
Since $\gamma_1$ is the highest noncompact root,
$|\gm_1|^2$ is the upper bound of curvature.
Let $\al=\sum m_j(\al)\al_j$,  from $(h_i, \al_j)=\fr 12 \delta_{ij}$ we have
\beq
\tta (e_{\al})=exp(2\pi i (\al, h_i))e_{\al}=(-1)^{m_i({\al})}e_\al
\eeq
then
\beq
\Delta_n=\{\al\in \Delta,  m_i(\al)\equiv 1 (mod 2)\}
\eeq
and $\gm_1$ be the highest root whose coefficient satisfies
$m_i\equiv 1 (mod 2)$.
This includes the cases of AIII, BI, DIII, CI, CII,  EI,  EII, EIII, EV, EVI, EVII,
EVIII, EIX, FI, FII, G and some cases of DI.
\par
(2)If $rank (\mfg)=rank (M)$,  then we can choose a Cartan subalgebra $\mfh$
whose real part satisfies $\mfh_0=\mfh_{\mfp}\subset \mfp$,
so the projection on $\mfh_{\mfp}$ of  any root is itself.
We choose the highest root $\al_0$,
$(\al_0, \al_0)$ is the curvature bound.
This is in the case of AI.
\par
(3) If $\tta=\tta_0$ is a special rotation,  $\tta_0(e_\al)=e_{\tta_0(\al)}$,
then $\Delta_n=0$,
$\mfh_\mfp$ is generated by
$$\{\al-\tta_0(\al),  \al\in \Delta_0^+ \}. $$
We can get the highest root which satisfies
$\al\in \Delta_0^+,  e_{\al}\in \mfp^\mbbc$,
the projection  is $\al_0'=\fr 12(\al-\tta_0(\al))$,
$(\al_0,  \al_0)$ is the upper bound of curvature.
This includes the cases of AII, EIV and some cases of DI.
\par
(4) For the other cases,  from the Cayley's transformation,
we can get $\{\gm_1, \gm_2, \cdots, \gm_s\}\cup \mfh_V$,
by calculating the projection on this subspace
we can give the upper bounds of sectional  curvature.
This includes some cases of DI. \par

Now  we calculate,   case by case,   the  upper bounds of sectional
curvature of irreducible symmetric spaces of compact type
according to the table in [4, p518]. \par

Let $\Pi=\{ \al_1, \al_2, \cdots, \al_l \}$ be a simple root system,
every root can be given as the integral linear combination
\beq
\al=m_1(\al)\al_1+m_2(\al)\al_2+\cdots+m_l(\al)\al_l.
\eeq
In the following calculation,  for four classical complex simple Lie algebras
 we imbed the root system into Euclidean spaces. \par

\setlength{\unitlength}{1cm}

(i)
The Dynkin diagram of $\mfa_l=\mf{sl}(l+1, \mbbc)$ is \\
\cl
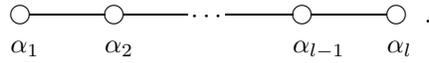
{
% Al
\begin{picture}(9, 1)
\put(0.125, 0.5){\circle{0.25}}
\put(0.25, 0.5){\line(1, 0){1}}
\put(1.375, 0.5){\circle{0.25}}
\put(1.5, 0.5){\line(1, 0){0.85}}
\put(2.4, 0.4){$\cdots$}
\put(2.85, 0.5){\line(1, 0){0.9}}
\put(3.875, 0.5){\circle{0.25}}
\put(4.0, 0.5){\line(1, 0){1.}}
\put(5.125, 0.5){\circle{0.25}}
\put(0,0){$\al_1$}
\put(1.25, 0){$\al_2$}
\put(3.75, 0){$\al_{l-1}$}
\put(5,0){$\al_l$}
\put(5.5,0.4){.}
\end{picture}
}

Let $\vep_j, 1\le j \le l+1$,  be an orthogonal base of $\mbbr^{l+1}$
with $|\vep_j|^2=\fr 1{2(l+1)}$.
The simple root system of $\mfa_l$ is
$$\Pi=\{ \al_j=\vep_j-\vep_{j+1},  j=1, 2, \cdots,  l \}, \q |\al_j|^2=\fr 1{l+1}. $$
The positive roots are $$\vep_i-\vep_j=\al_i+\cdots+\al_{j-1},  i<j. $$
The highest root is $$\delta=\vep_1-\vep_{l+1}=\al_1+\al_2+\cdots+\al_l. $$
\par

For AI,  M=$SU(n)/SO(n),  \mfg=\mfa_l=\mf{sl}(n, \mbbc),  l=n-1,
rank (M)=rank (\mfg)=l,  \theta(Z)=-Z^T,  Z\in \mfg$.
If n is even,  then $n=2m, \tta=\tta_0 exp(2\pi i h_m)$,
if n is odd,  then  $n=2m+1,  \tta=\tta_0$,  where $\tta_0$ is the special rotation.
In this case ,  let  $\mfh_0$ consists of all real diagonal matrixes of trace 0,
then $\mfh_0$ is a  Cartan subalgebra of $\mfg_0$,  as well as a
reduced Cartan subalgebra of $(\mfg_0,  \tta)$,  $\mfh_0\subset \mfp$,
so the length of restricted root $\al'$ is same as the length of $\al$.
Since for $\mfa_l$,  all simple roots (as well as all roots) have the same length,  $(\al,  \al)=\fr 1{l+1}$,
We see that  the upper bound of curvature is $\fr 1n$. \par

For AII,  $M=SU(2n)/Sp(n)$, $\mfg=\mfa_l,  l=2n-1,  \theta=\theta_0,  \tta_0(Z)=-J_nZ^TJ_n^{-1}$ is
special rotation,
$$\tta(\al_j)=\al_{l+1-j}=\al_{2n-j},  j=1, 2, \cdots,  2n-1.$$
If $$\al=m_1\al_1+\cdots + m_n\al_n+\cdots+ m_{2n-1}\al_{2n-1}$$
then
\beqs
&&\tta(\al)=m_{2n-1}\al_1+\cdots + m_n\al_n+\cdots+ m_1\al_{2n-1}\\
&&\tta(e_\al)=e_{\tta(\al)}
\eeqs
In this case $\Delta_n$ is empty set,  the reduced Cartan subalgebra $\mfh_{\mfp}$ is generated by
$\{ \al-\theta (\al),  \al\in \Delta \}$,
the  base is
$$\{ \sigma_i=\al_i-\tta(\al_i)=\vep_i-\vep_{i+1}-\vep_{2n-i}+\vep_{2n-i+1},  1\le i\le n-1 \}. $$
$\tta_0(\al)=\al$ if and only if $m_i(\al)=m_{2n-i}(\al),  i=1, 2, \cdots,  n-1$.
The positive root is in the form of $\al=\vep_i-\vep_j=\al_i+\cdots+\al_{j-1}$.
The highest root which satisfies $\tta_0(\al)\ne \al$ is
$$\al_0=\vep_1-\vep_{2n-1}=\al_1+\cdots+\al_{2n-2}. $$
We get the orthogonal  base $\{ \tilde\sigma_i \}$ from $\{\sigma_i\}$ such that
\beqs
& \tilde\sigma_1=\sigma_1\\
&\tilde\sigma_2=2\sigma_2+\sigma_1=\vep_1+\vep_2-2\vep_3-2\vep_{2n-2}+\vep_{2n-1}+\vep_{2n}\\
&\cdots\\
&(\tilde\sigma_i,  \al_0)=0,  i>1.
\eeqs
The projection on $\{ \tilde\sigma_i,  1\le i\le n-1\}$ of $\al_0$ is
$$\fr{(\al_0,  \tilde\sigma_1)}{(\tilde \sigma_1,  \tilde\sigma_1)}\tilde\sigma_1$$
with the squared length
$$\fr{(\al_0,  \tilde\sigma_1)^2}{(\tilde \sigma_1,  \tilde\sigma_1)}=\fr 12 (\al_0, \al_0) =\fr 1{4n}. $$
We see that the upper bound of curvature is $\fr 1{4n}$. \par

For $AIII$, $M=SU(p+q)/S(U_p\tm U_q),  \mfg=\mfa_l,  l=p+q-1,  rank (M)=min(p, q),
\tta(Z)=-I_{p, q}X^TI_{p, q},
\tta=exp(2\pi i ad h_i),  rank(\mfk)=rank(\mfg),
\al_0=\al_1+\cdots+\al_l\in \Delta_n$
is the highest root.
In  the Cayley's  transformation we have $\gm_1=\al_0$,
but $|\al_0|^2=\fr 1{l+1}$ attains the maximal value,
so the upper bound of curvature is $|\al_0|^2=\fr 1{p+q}$. \par

For BDI, $M=SO(p+q)/SO(p)\tm SO(q), \tta(Z)=I_{p, q}Z^TI_{p, q}$,
these are cases of BI and DI according to $p+q$ being a odd or a even. \par

(ii)
The Dynkin diagram of $\mfb_l=\mf{so}(2l+1, \mbbc)$ is \\
\cl
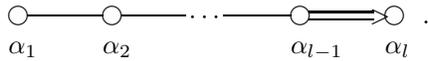
{
% Bl
\begin{picture}(9, 1)
\put(0.125, 0.5){\circle{0.25}}
\put(0.25, 0.5){\line(1, 0){1}}
\put(1.375, 0.5){\circle{0.25}}
\put(1.5, 0.5){\line(1, 0){0.85}}
\put(2.4, 0.4){$\cdots$}
\put(2.85, 0.5){\line(1, 0){0.9}}
\put(3.875, 0.5){\circle{0.25}}
\put(4.0, 0.55){\line(1, 0){0.85}}
\put(4.0, 0.45){\line(1, 0){0.85}}
\put(4.8, 0.4){$>$}
\put(5.125, 0.5){\circle{0.25}}
\put(0, 0){$\al_1$}
\put(1.25, 0){$\al_2$}
\put(3.75, 0){$\al_{l-1}$}
\put(5, 0){$\al_l$}
\put(5.5 ,0.4){.}
\end{picture}
}\par
Let $\vep_j, 1\le j \le l$,  be an orthogonal base of $\mbbr^{l},  |\vep_j|^2=\fr 1{2(2l-1)}$,
the simple root system of $\mfb_l$ is
$$\Pi=\{ \al_j=\vep_j-\vep_{j+1},  j=1, 2, \cdots,  l-1,  \al_l=\vep_l \}. $$
The positive roots are $$\vep_i,  \vep_i\pm \vep_j,  i<j$$
where
\beqs
&\vep_i=\al_i+\cdots+\al_l    \\
&\vep_i-\vep_j=\al_i+\cdots+\al_j \\
&\vep_i+\vep_j=\al_i+\cdots+\al_{j-1}+2(\al_j+\cdots+\al_l).
\eeqs
The highest root is
$$\delta=\vep_1+\vep_2=\al_1+2\al_2+\cdots+2\al_l. $$
For BI,  $M=SO(p+q)/SO(p)\tm SO(q), p+q=2l+1$, $\mfg=\mfb_l, rank(\mfk)=rank(\mfg), \\
\tta=exp(2\pi i ad h_i)$ is a inner automorphism
$$\Delta_n=\{\al\in \Delta,  m_i(\al) \equiv  1 (mod 2) \}. $$
The highest noncompact root $\al_0$ is the highest root whose coefficient $m_i$ is a odd.
We get it in these two cases
\beqs
& \vep_1+\vep_{i+1}=\al_1+\cdots+\al_i+2(\al_{i+1}+\cdots+\al_l), \q i<l,\\
& \vep_1=\al_1+\al_2+\cdots+\al_l.
\eeqs
So the upper bound of curvature is
$\fr 1{2l-1}=\fr 1{p+q-2}(i< l),  \fr 1{2(2l-1)}=\fr 1{2(p+q-2)}(i=l)$.
The later case corresponds  to $SO(2n+1)/SO(2n)$ which is of rank 1.\par

(iii) The Dynkin diagram of  $\mfd_l=\mf{so}(2l, \mbbc)$ is \\
\cl
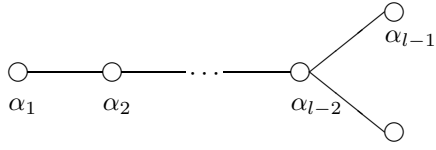
{
% dl
\begin{picture}(9, 2)
\put(0.125, 0.5){\circle{0.25}}
\put(0.25, 0.5){\line(1, 0){1}}
\put(1.375, 0.5){\circle{0.25}}
\put(1.5, 0.5){\line(1, 0){0.85}}
\put(2.4, 0.4){$\cdots$}
\put(2.85, 0.5){\line(1, 0){0.9}}
\put(3.875, 0.5){\circle{0.25}}
\put(4.0, 0.5){\line(5, 4){1}}
\put(4.0, 0.5){\line(5, -4){1}}
\put(5.125, 1.3){\circle{0.25}}
\put(5.125, -0.3){\circle{0.25}}
\put(5., 0.9){$\al_{l-1}$}
\put(5., -0.7){$\al_{l}$}
\put(0, 0){$\al_1$}
\put(1.25, 0){$\al_2$}
\put(3.75, 0){$\al_{l-2}$}
%\put(9.  , 0.2){\Large{$\mfd_l$}}
\end{picture}
}
\par
\vspace{0.5cm}
Let $\vep_j, 1\le j \le l$,  be an orthogonal base of $\mbbr^{l}, |\vep_j|^2=\fr 1{4(l-1)}$.
The simple root system is
$$\Pi=\{ \al_j=\vep_j-\vep_{j+1},  j=1, 2, \cdots,  l-1,  \al_l=\vep_{l-1}+\vep_l \}. $$
The positive roots are
$$\vep_i\pm \vep_j,  i<j$$
where
\beqs
&\vep_i-\vep_j&=\al_i+\cdots+\al_{j-1}\\
&\vep_i+\vep_j&=\al_i+\cdots+\al_{l-2}+\al_{j}+\cdots+\al_l\\
&    &=\al_i+\cdots+\al_{j-1}+2(\al_{j}+\cdots+\al_{l-2})+\al_{l-1}+\al_l.
\eeqs
The highest root is
$$\delta=\vep_1+\vep_{2}=\al_1+2\al_2+\cdots+2\al_{l-2}+\al_{l-1}+\al_l. $$

For DI,  $M=SO(p+q)/SO(p)\tm SO(q), p+q=2l$, $\mfg=\mfd_l$. \par
$\tta$ can be a inner automorphism or a outer automorphism. \par
If $\tta=exp(2\pi i ad h_i), 1\le i\le [\fr l2]$,
we can get the highest root which satisfies $m_i(\al)=1$,
$$
\al_0=\vep_1+\vep_{i+1}=\al_1+\cdots+\al_i+2(\al_{i+1}+\cdots+\al_{l-2})+\al_{l-1}+\al_l.
$$
So the maximum of curvature is $2\cdot \fr 1{4(l-1)}=\fr 1{p+q-2}$. \par

If $\theta=\theta_0$ is a outer involution,
$$\tta_0(\al_i)=\al_i,  1\le i\le l-2,  \tta_0(\al_{l-1})=\al_l, \tta_0(\al_l)=\al_{l-1}$$
then $\Delta_n$ is empty,  the reduced Cartan subalgebra $\mfh_\mfp$ is generated by
$$\{ \sigma_1=\al_{l-1}-\theta (\al_{l-1})=\al_{l-1}-\al_l\}. $$
$rank(M)=dim \mfh_\mfp=1, \tta_0(\al)=\al$ if and only if $m_{l-1}(\al)=m_l(\al)$. \\
The highest root which satisfies $m_{l-1}(\al)\ne m_l(\al)$ is
$$\al_0=\vep_1-\vep_l=\al_1+\cdots+\al_{l-1}. $$
The projection of $\al_0$  is
$$\fr{(\al_0, \sigma_1)}{(\sigma_1, \sigma_1)}\sigma_1$$
with the squared length
$\fr{(\al_0, \sigma_1)^2}{(\sigma_1, \sigma_1)}=\fr 12 (\al_0, \al_0)=\fr 1{2(p+q-2)}$.
We see that the maximum of curvature is $\fr 1{2(p+q-2)}$. \par

If $\theta=\theta_0exp(2\pi i ad h_i),  1\le i\le [\fr l2]$,
from Cayley's transformation,  we know  that there exist
$\gm_1, \gm_2, \cdots,  \gm_{r-1}$,
the reduced Cartan subalgebra $\mfh_\mfp$ is generated by
$$\{ \gm_i,  i=1, 2, \cdots,  r-1,   \q \al_{l-1}-\al_l \}. $$
The projection of root $\gm_1$ on $\mfh_\mfp$ is itself,
but $|\gm_1|$ attains the maximum of length of roots,
we see that the maximum of curvature   is $\fr 1{p+q-2}$. \par

For DIII,  $M=SO(2n)/U(n),  \mfg=\mfg_l,  l=n,  \tta=exp(2\pi i ad h_i),  i=l$ is inner involution,
the highest noncompact root is
$$ \vep_1+\vep_2=\al_1+2(\al_{2}+\cdots+\al_{l-1})+\al_{l-1}+\al_l. $$
We see  that  the maximum of curvature  is $\fr 1{2(l-1)}=\fr 1{2n-2}$. \par

(iv) The Dynkin diagram of $\mfc_l=\mf{sp}(l, \mbbc)$ is \\
\cl{
% cl
\begin{picture}(9, 1)
\put(0.125, 0.5){\circle{0.25}}
\put(0.25, 0.5){\line(1, 0){1}}
\put(1.375, 0.5){\circle{0.25}}
\put(1.5, 0.5){\line(1, 0){0.85}}
\put(2.4, 0.4){$\cdots$}
\put(2.85, 0.5){\line(1, 0){0.9}}
\put(3.875, 0.5){\circle{0.25}}
\put(4, 0.4){$<$}
\put(4.15, 0.55){\line(1, 0){0.85}}
\put(4.15, 0.45){\line(1, 0){0.85}}
\put(5.125, 0.5){\circle{0.25}}
\put(0, 0){$\al_1$}
\put(1.25, 0){$\al_2$}
\put(3.75, 0){$\al_{l-1}$}
\put(5, 0){$\al_l$}
\put(5.5, 0.4){.}
\end{picture}
}
\par
Let $\vep_j, 1\le j \le l$,  be an orthogonal base of
$\mbbr^{l}, |\vep_j|^2=\fr 1{4(l+1)}$,
the simple root system is
$$\Pi=\{ \al_j=\vep_j-\vep_{j+1},  j=1, 2, \cdots,  l-1,  \al_l=2\vep_l \}. $$
The positive roots are
$$2\vep_i,  \vep_i\pm \vep_j,  i<j$$
where
\beqs
&2\vep_i  = & 2(\al_i+\cdots+\al_{l-1})+\al_l\\
&\vep_i-\vep_j  = & \al_i+\cdots+\al_{j-1}\\
&\vep_i+\vep_j =& \al_i+\cdots+\al_{j-2}+2(\al_{j}+\cdots+\al_{l-1})+\al_l.
\eeqs
The highest root is
$$\delta=2\vep_1=2\al_1+2\al_2+\cdots+2\al_{l-1}+\al_l. $$
\par
For CI or CII,  $\mfg=\mfc_l, \tta=exp(2\pi i ad h_i),  1\le i\le [\fr {l-1}{2}]+1$ or $i=l$. \par
For CI,  $i=l$,  the noncompact root satisfies  $m_l(\al)$ is odd,
we see that the highest noncompact root is
$  \al_0=2\vep_1=2(\al_1+\al_2+\cdots +\al_{l-1})+\al_l. $
The maximum of curvature  is  $(\al_0,  \al_0)=\fr 1{n+1}$\par
For CII,  $l=p+q,  i<l$,  the noncompact root satisfies  $m_i(\al)$ is odd,
we see that the highest noncompact root is
$\al_0=\vep_i+\vep_{i+1}=\al_i+2(\al_{i+1}+\cdots+\al_{l-1})+\al_l. $
The maximum of curvature  is
$(\al_0, \al_0)=\fr 1{2(l+1)}=\fr 1{2(p+q+1)}. $\par

Now we  consider the cases of exception type.  \par
(v)For type E,
$\mfg=\mfe_6, \mfe_7, \mfe_8$,
we draw the corresponding Dynkin diagram as follows. \par
% E6
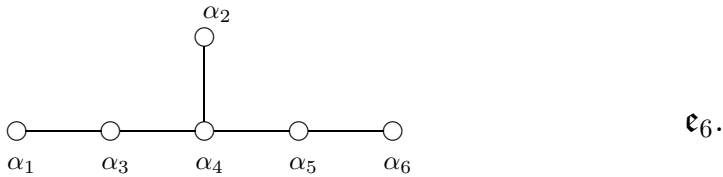
\begin{picture}(9, 2)
\put(0.125, 0.5){\circle{0.25}}
\put(0.25, 0.5){\line(1, 0){1}}
\put(1.375, 0.5){\circle{0.25}}
\put(1.5, 0.5){\line(1, 0){1.}}
\put(2.625, 0.5){\circle{0.25}}
\put(2.75, 0.5){\line(1, 0){1.}}
\put(3.875, 0.5){\circle{0.25}}
\put(4.0, 0.5){\line(1, 0){1.}}
\put(5.125, 0.5){\circle{0.25}}
\put(0, 0){$\al_1$}
\put(1.25, 0){$\al_3$}
\put(2.5, 0){$\al_4$}
\put(3.75, 0){$\al_5$}
\put(5, 0){$\al_6$}
\put(2.5, 0.625){ \line(0, 1){1} }
\put(2.625, 1.75){\circle{0.25}}
\put(2.6, 2){$\al_2$}
\put(9, 0.5){\Large{$\mfe_6.$}}
\end{picture}
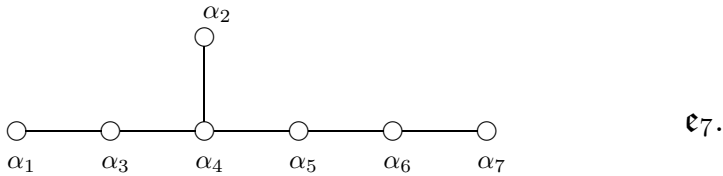
\par
\begin{picture}(9, 2.5)
\put(0.125, 0.5){\circle{0.25}}
\put(0.25, 0.5){\line(1, 0){1}}
\put(1.375, 0.5){\circle{0.25}}
\put(1.5, 0.5){\line(1, 0){1.}}
\put(2.625, 0.5){\circle{0.25}}
\put(2.75, 0.5){\line(1, 0){1.}}
\put(3.875, 0.5){\circle{0.25}}
\put(4.0, 0.5){\line(1, 0){1.}}
\put(5.125, 0.5){\circle{0.25}}
\put(5.25, 0.5){\line(1, 0){1.}}
\put(6.375, 0.5){\circle{0.25}}
\put(0, 0){$\al_1$}
\put(1.25, 0){$\al_3$}
\put(2.5, 0){$\al_4$}
\put(3.75, 0){$\al_5$}
\put(5, 0){$\al_6$}
\put(6.25, 0){$\al_7$}
\put(2.5, 0.625){ \line(0, 1){1} }
\put(2.625, 1.75){\circle{0.25}}
\put(2.6, 2){$\al_2$}
\put(9, 0.5){\Large{$\mfe_7.$}}
\end{picture}
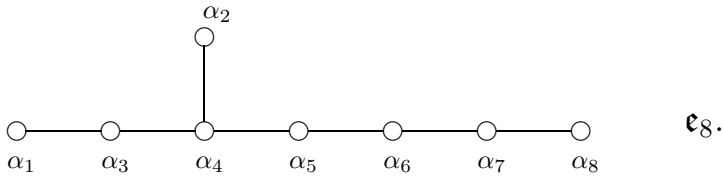
\par
% E8
\begin{picture}(9, 2)
\put(0.125, 0.5){\circle{0.25}}
\put(0.25, 0.5){\line(1, 0){1}}
\put(1.375, 0.5){\circle{0.25}}
\put(1.5, 0.5){\line(1, 0){1.}}
\put(2.625, 0.5){\circle{0.25}}
\put(2.75, 0.5){\line(1, 0){1.}}
\put(3.875, 0.5){\circle{0.25}}
\put(4.0, 0.5){\line(1, 0){1.}}
\put(5.125, 0.5){\circle{0.25}}
\put(5.25, 0.5){\line(1, 0){1.}}
\put(6.375, 0.5){\circle{0.25}}
\put(6.5, 0.5){\line(1, 0){1.}}
\put(7.625, 0.5){\circle{0.25}}
\put(0, 0){$\al_1$}
\put(1.25, 0){$\al_3$}
\put(2.5, 0){$\al_4$}
\put(3.75, 0){$\al_5$}
\put(5, 0){$\al_6$}
\put(6.25, 0){$\al_7$}
\put(7.5, 0){$\al_8$}
\put(2.5, 0.625){ \line(0, 1){1} }
\put(2.625, 1.75){\circle{0.25}}
\put(2.6, 2){$\al_2$}
\put(9, 0.5){\Large{$\mfe_8.$}}
\end{picture}
\par
In these cases all the roots have the same length.
All involutions are inner automorphism except for $\mfe_6$.
From formula (\ref{ci}),
we get the squared length of roots are  $\fr 1{12},  \fr 1{18}, \fr 1{30}$, respectively.

\par
We list the positive roots of $\mfe_6$ as follows:
\beqs
&&\al_1, \q \al_2, \q \al_3, \q \al_4, \q \al_5, \q \al_6\\
&&\al_1+\al_3, \q \al_3+\al_4, \q \al_2+\al_4, \q\al_4+\al_5, \q\al_5+\al_6\\
&& \al_1+\al_3+\al_4, \q \al_3+\al_4+\al_2, \q\al_3+\al_4+\al_5\\
    &&\al_4+\al_5+\al_2, \q\al_4+\al_5+\al_6\\
&& \al_1+\al_3+\al_4+\al_2, \q \al_1+\al_3+\al_4+\al_5, \q\al_2+\al_3+\al_4+\al_5\\
  &&\al_3+\al_4+\al_5+\al_6, \q\al_2+\al_4+\al_5+\al_6\\
&& \al_1+\al_2+\al_3+\al_4+\al_5, \q \al_1+\al_3+\al_4+\al_5+\al_6\\
    &&\al_2+\al_3+\al_4+\al_5+\al_6, \q\al_2+\al_3+2\al_4+\al_5\\
&& \al_1+\al_2+\al_3+\al_4+\al_5+\al_6, \q \al_1+\al_2+\al_3+2\al_4+\al_5\\
    &&\al_2+\al_3+2\al_4+\al_5+\al_6\\
&& \al_1+\al_2+\al_3+2\al_4+\al_5+\al_6, \q \al_1+\al_2+2\al_3+2\al_4+\al_5\\
    &&\al_2+\al_3+2\al_4+2\al_5+\al_6\\
&& \al_1+\al_2+2\al_3+2\al_4+\al_5+\al_6, \q \al_1+\al_2+\al_3+2\al_4+2\al_5+\al_6\\
&& \al_1+\al_2+2\al_3+2\al_4+2\al_5+\al_6\\
&& \al_1+\al_2+2\al_3+3\al_4+2\al_5+\al_6\\
&& \al_1+2\al_2+2\al_3+3\al_4+2\al_5+\al_6.
\eeqs
For EIV,  $\mfg=\mfe_6,  \theta=\tta_0$ is outer automorphism
\beq
& \tta_0(\al_1)=\al_6, \tta_0(\al_2)=\al_2, \tta_0(\al_3)=\al_5 \\
& \tta_0(\al_4)=\al_4,  \tta_0(\al_5)=\al_3, \tta_0(\al_6)=\al_1.
\eeq
The reduced Cartan subalgebra $\mfh_\mfp$ is generated by
$$\{ \al_1-\al_6,  \al_3-\al_5 \}. $$
$\tta_0(\al)=\al$ if and only if $m_1(\al)=m_6(\al),  m_3(\al)=m_5(\al). $
The highest root which satisfies $\tta_0(\al)\ne \al$ is
$$\al_0=\al_1+\al_2+2\al_3+2\al_4+\al_5+\al_6.$$
A orthogonal base of $\mfh_\mfp$ is
$$\{ \sigma_1=\al_1-\al_6, \sigma_2=\al_3-\al_5+\fr 12(\al_1-\al_6) \}. $$
The projection on $\mfh_\mfp$ of $\al_0$ is
$$\fr{(\al_0,  \sigma_1)}{(\sigma_1, \sigma_1)}\sigma_1+\fr{(\al_0,  \sigma_2)}{(\sigma_2, \sigma_2)}\sigma_2.$$
We see that  the maximum of curvature  is
$$\fr{(\al_0,  \sigma_1)^2}{(\sigma_1, \sigma_1)}+\fr{(\al_0,  \sigma_2)^2}{(\sigma_2, \sigma_2)}
=\fr 12{(\al_0, \al_0)}=\fr 1{24}. $$\par

For EI,  EII, EIII,  we need the Cayley's transformation,  there exists at least one
$\gm_1$,  but all roots have the same length,  $|\gm_1|^2=\fr 1{12}$
attains the maximum.
We see that  the maximum of curvature  is $|\gm_1|^2=\fr 1{12}$.
\par
For EV, EVI, EVII, $\mfg=\mfe_7$,  all involutions are inner automorphisms,
with the same reason as previous case we get
the maximum of curvature  is $\fr 1{18}$.
\par
For EVIII, EIX, $\mfg=\mfe_8$,  all involutions are inner automorphisms,
with the same reason as previous case we get
the maximum of curvature  is $\fr 1{30}$.
\par
(vi)
The Dynkin diagram of $\mff_4$ is \\
\cl
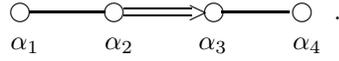
{
% f4
\begin{picture}(9, 1)
\put(0.125, 0.5){\circle{0.25}}
\put(0.25, 0.5){\line(1, 0){1}}
\put(1.375, 0.5){\circle{0.25}}
\put(1.5, 0.55){\line(1, 0){0.9}}
%\put(1.5, 0.5){\line(1, 0){0.9}}
\put(1.5, 0.45){\line(1, 0){0.9}}
\put(2.35, 0.40){$>$}
\put(2.70, 0.5){\circle{0.25}}
\put(2.8, 0.5){\line(1, 0){0.9}}
\put(3.875, 0.5){\circle{0.25}}
\put(0, 0){$\al_1$}
\put(1.25, 0){$\al_2$}
\put(2.5, 0){$\al_3$}
\put(3.75, 0){$\al_4$}
\put(4.3, 0.4){.}
\end{picture}
}
For FI, FII, $\mfg=\mff_4,  |\al_1|^2=|\al_2|^2=\fr 1{9},  |\al_3|^2=|\al_4|^2=\fr 1{18}. $
We list all positive roots
\beqs
&& \al_1, \q \al_2, \q \al_3, \q \al_4\\
&& \al_1+\al_2, \q \al_2+\al_3, \q \al_3+\al_4\\
&& 2\al_2+\al_3, \q \al_1+\al_2+\al_3, \q \al_2+\al_3+\al_4\\
&& \al_1+\al_2+2\al_3, \q \al_1+\al_2+\al_3+\al_4, \q \al_2+2\al_3+\al_4\\
&& \al_1+2\al_2+2\al_3, \q \al_2+2\al_3+2\al_4, \q \al_1+\al_2+2\al_3+\al_4\\
&& \al_1+\al_2+2\al_3+2\al_4, \q \al_1+2\al_2+2\al_3+\al_4\\
&& \al_1+2\al_2+2\al_3+2\al_4, \q \al_1+2\al_2+3\al_3+\al_4\\
&& \al_1+2\al_2+3\al_3+2\al_4\\
&& \al_1+2\al_2+4\al_3+2\al_4\\
&& \al_1+3\al_2+4\al_3+2\al_4\\
&& 2\al_1+3\al_2+4\al_3+2\al_4.
\eeqs
All involutions are inner automorphisms,  $\tta=exp(2\pi i ad h_i),  i=1, 4.$\\
For $i=1$,  the highest noncompact root which satisfies  $m_1=1 (mod \ 2)$ is
$$\al_0=\al_1+3\al_2+4\al_3+2\al_4,  |\al_0|^2 =\fr 1{9}. $$
For $i=4$,  the highest noncompact root which satisfies $m_4=1 (mod \ 2)$ is
$$\al_0=\al_1+2\al_2+3\al_3+\al_4,  |\al_0|^2 =\fr 1{18}. $$
We get the maximum of curvature is $\fr 19$ for FI,  $\fr 1{18}$ for  FII. \par

(vii)
The Dynkin diagram of $\mfg_2$ is \\
\cl
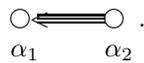
{
% g2
\begin{picture}(9, 1)
\put(0.125, 0.5){\circle{0.25}}
\put(0.25, 0.40){$<$}
\put(0.35, 0.5){\line(1, 0){0.9}}
\put(0.35, 0.55){\line(1, 0){0.9}}
\put(0.35, 0.45){\line(1, 0){0.9}}
\put(1.375, 0.5){\circle{0.25}}
\put(0, 0){$\al_1$}
\put(1.25, 0){$\al_2$}
\put(1.7,0.4){.}
\end{picture}
}
For G, $\mfg=\mfg_2,  |\al_1|^2=\fr 1{12},  |\al_2|^2=\fr 1{4}$,
$\tta=exp(2\pi i ad h_i),  i=2$. We list all positive roots
\beqs
&& \al_1, \q \al_2\\
&& \al_1+\al_2, \q  2\al_1+\al_2\\
&& 3\al_1+\al_2, \q  3\al_1+2\al_2.
\eeqs
The highest noncompact root which satisfies $m_2=1 (mod \ 2)$ is
 $$\al_0=3\al_1+\al_2, |\al_0|^2=\fr 14.$$
We get the maximum of curvature is $\fr 14$. \par

Now we get the upper bounds of sectional curvature for all cases of irreducible Riemannian
symmetric spaces of compact type. In summary,  we have the following table.
\par
\vspace{0.3cm}
\cl{Table 5.1}
\vspace{0.3cm}
\cl{
\begin{tabular}{|c|c|c|c|c|c|c|}
  \hline
 Type &   compact type & rank & dimension & bound\\
 \hline
 $AI$ &   $SU(n)/SO(n)$  & $n-1$  & $\fr 12(n-1)(n+2)$  &  $\fr 1n$\\
 \hline
 $AII$ &  $SU(2n)/Sp(n)$  & $n-1$  & $(n-1)(2n+1)$  &  $\fr 1{4n}$ \\
 \hline
 $AIII$ & $SU(p+q))/S(U_p\tm U_q))$  & $min(p, q)$  & $2pq$  &  $\fr 1{p+q}$ \\
 \hline
 $BDI$ & $SO(p+q))/SO(p)\tm SO(q)$ & $min(p, q)>1$  & $pq$  &   $\fr 1{p+q-2}$\\
   &  & \q\q\q\q=1  &   &  $\fr 1{2(p+q-2)}$\\
 \hline
 $DIII$ & $SO(2n))/U(n)$ & $[\fr 12n ] $  & $n(n-1)$  &  $\fr 1{2n-2}$ \\
 \hline
 $CI$ & $Sp(n)/U(n)$ & $n$  & $n(n+1)$  &  $\fr 1{n+1}$ \\
 \hline
 $CII$ & $Sp(p+q))/Sp(p)\tm Sp(q)$ & $min(p, q)$  & $4pq$  &  $\fr 1{2(p+q+1)}$ \\
 \hline
 $EI-III$ & $$ & $6/4/2$  & $42/40/32$  &  $\fr 1{12}$ \\
 \hline
 $EIV$ & $$ & $2$  & $26$  &  $\fr 1{24}$ \\
 \hline
 $EV-VII$ & $$ & $7/4/3$  & $70/64/54$  &  $\fr 1{18}$ \\
 \hline
 $EVIII-IX$ & $$ & $8/4$  & $128/112$  &  $\fr 1{30}$ \\
 \hline
 $FI$ & $$ & $4$  & $28$  &  $\fr 19$ \\
 \hline
 $FII$ & $$ & $1$  & $16$  &  $\fr 1{18}$ \\
 \hline
 $G$ & $$ & $2$  & $8$  &  $\fr 14$ \\
 \hline
 \end{tabular}
 }
\par
\remark  {Finally,  we would like to emphasize the normality for
 a compact irreducible symmetric space $M=U/K$  of compact type.
   The  invariant Riemannian inner product at $o=\pi(e)$ is \beq \qq
<X,Y>=-c\cdot B_{\mfu}(X,Y)=-c \cdot B_{\mfg}(X,Y),\q X,Y\in \mfu,
\mfg=\mfu^\mbbc, \eeq where $c>0$ is a positive constant,
$B_{\mfu}(X,Y)$ is the Killing form of $\mfu$,  and it's the
restriction of the Killing form of $\mfg$. In above calculations,
we choose c=1,  see (\ref{innp}). Under this normalization, the
Ricci curvature is $\fr 12$ for any case. }

{\bf Acknowledgements.}
\thanks
{ I am  grateful to my supervisor, Prof. Y.L.Xin, for his constant
guidance. He raised this problem for my investigation. Prof. Xin
also gave me many valuable suggestions  and pointed out some
mistakes in the course  of this work. }

\end{document}